\input amssym.def
\input psfig
\input epsf

\let \blskip = \baselineskip
\parskip=1.2ex plus .2ex minus .1ex

\tabskip 20pt
\tolerance = 1000
\pretolerance = 50
\newcount\itemnum
\itemnum = 0
\overfullrule = 0pt

\def\title#1{\bigskip\centerline{\bigbigbf#1}}
\def\author#1{\bigskip\centerline{\bf #1}\smallskip}
\def\address#1{\centerline{\it#1}}
\def\abstract#1{\vskip1truecm{\narrower\noindent{\bf Abstract.} #1\bigskip}}

\def\sp{\bigskip}
\def\nosp{\vskip -\the\blskip plus 1pt minus 1pt}
\def\upsp{\nosp\medskip}
\def\br{\hfil\break} 
\def\ti{\br \hglue \the \parindent}

\def\ce#1{\LP\centerline{#1}}

\def\skipit#1{}
\def\mdag{\raise 3pt\hbox{\dag}}

\def\XP{\par\noindent\hang}
\def\LP{\par\noindent}
\def\BP[#1]{\par\item{[#1]}}
\def\SH#1{\sp\vskip\parskip\leftline{\bigbf #1}\nobreak}

\def\TH#1{\sp\XP{\bf THEOREM\ \shead#1}}
\def\LM#1{\sp\XP{\bf LEMMA\ \shead#1}}
\def\DF#1{\sp\XP{\bf DEFINITION\ \shead#1}}
\def\PR#1{\sp\XP{\bf PROPOSITION\ \shead#1}}
\def\CO#1{\sp\XP{\bf COROLLARY\ \shead#1}}

\def\PF{\LP{\bf Proof:\ }}
\def\NX{\advance\itemnum by 1 \sp\LP {\bf \shead \the\itemnum.\ }}
\def\qed{\null\nobreak\hfill\hbox{${\vrule width 5pt height 6pt}$}\par\sp}

\def\cart{\>\hbox{${\vcenter{\vbox{
    \hrule height 0.4pt\hbox{\vrule width 0.4pt height 4.5pt
    \kern4pt\vrule width 0.4pt}\hrule height 0.4pt}}}$}\>}
\def\bxmu{\>\hbox{${\vcenter{\vbox {
    \hrule height 0.4pt\hbox{\vrule width 0.4pt height 4pt
    \hskip -1.3pt\lower 1.8pt\hbox{$\times$}\negthinspace\vrule width 0.4pt}
    \hrule height 0.4pt}}}$}\>}

\def\lin#1{\hbox to #1true in{\hrulefill}}


\def\foot#1{\raise 6pt \hbox{#1} \kern -3pt}

\def\fig #1 #2 #3 #4 #5 {\sp \ce{ {\epsfbox[#1 #2 #3 #4]{figs/#5.ps}} }}

\def\gpic#1{#1 \sp\ce{\box\graph} \medskip} 


\def\JGT{{\it J.\ Graph Theory}}

\def\DM{{\it Discrete Math.{}}}
\def\DAM{{\it Discrete Appl.\ Math.{}}}

\def\SIAP{{\it SIAM J.\ Appl.\ Math.{}}}

\def\SIREV{{\it SIAM Review{}}}

\def\al{\alpha}



\def\nul{\hbox{\O}}    
		
\def\esub{\subseteq}

\def\comp{\circ}

 \def\Imp{\Rightarrow} 
  

\def\({\left(}	\def\){\right)}


\def\UE#1#2#3{\bigcup_{#1 = #2} ^ {#3}}

\def\VEC#1#2#3{#1_{#2},\ldots,#1_{#3}}

\def\st{\colon\;} 

\def\SET#1:#2{\{#1\colon\;#2\}}


		






\magnification=\magstep1
\vsize=9.0 true in
\hsize=6.5 true in
\headline={\hfil\ifnum\pageno=1\else\folio\fi\hfil}
\footline={\hfil\ifnum\pageno=1\folio\else\fi\hfil}

\parindent=20pt
\baselineskip=12pt
\parskip=.5ex  

\def\shead{ }

\font\bigbf = cmb10 scaled \magstep1

\font\bigbigbf = cmb10 scaled \magstep2



\def\ba#1{\beta(\alpha(#1))}
\title{LINE DIGRAPHS AND COREFLEXIVE VERTEX SETS}
\author{Xinming Liu}
\address{Department of Industrial Economics and Systems Engineering}
\address{Tianjin University, Tianjin 300072, PR China, xliu@writeme.com}
\author{Douglas B. West}
\address{Mathematics Department}
\address{University of Illinois, Urbana, IL 61801-2975, {west@math.uiuc.edu}}
\vfootnote{}{\br
   Running head: LINE DIGRAPHS AND COREFLEXIVE SETS \br
}
\abstract{
The concept of coreflexive set is introduced to study the structure of
digraphs.  New characterizations of line digraphs and $n$th-order line
digraphs are given.  Coreflexive sets also lead to another natural way of
forming an intersection digraph from a given digraph.}

\SH
{1. INTRODUCTION}
In this paper we introduce new structural concepts about directed graphs.  We
use $V(D)$ and $E(D)$ for the vertex set and edge set of a digraph $D$.  We
allow loops but no multiple edges, unless specified otherwise.  We write $(u,v)$
or $uv$ for an edge from $u$ to $v$, with {\it tail} $u$ and {\it head} $v$.

Let $D$ be a digraph with multiple edges allowed.  The {\it line digraph} of $D$
is a digraph $L(D)$ (without multiple edges) such that $V(L(D))=E(D)$, and
$L(D)$ has an edge from $e$ to $f$ if and only if the head of $e$ is the tail of
$f$.  Introduced in [4], line digraphs were studied also in [2,5,6,7] ([5]
contains a survey).  Beineke and Zamfirescu [1] studied the $n$th-order line
digraph $L^n(D)$ (obtained by iterating the line digraph operation), and they
characterized the second-order line digraphs.  We obtain several
characterizations of line digraphs and a characterization of $n$th-order line
digraphs.

To study the structure of digraphs, we introduce the concept of ``coreflexive
set''.  We use $\al(v)$ to denote the set of successors of a vertex $v$, also
commonly written as $N^+(v)$ (the ``out-neighbor'' set).  We use $\beta(v)$ to
denote the set of predecessors of $v$, also commonly written
as $N^-(v)$ (the ``in-neighbor'' set).  Our use of $\al$ and $\beta$ suggests
the words ``after'' and ``before''.  We extend this notation to sets of
vertices:  $\al(S)$ is the set of vertices to which there is an edge from
at least one vertex of $S$, and $\beta(S)$ is the set of vertices from which
there is an edge to at least one vertex of $S$.

\DF{}
A {\it coreflexive set} in a digraph (henceforth {\it coreset})
is either 1) the collection of all sinks
in the digraph, or 2) a minimal nonempty set $U$ such that $U=\ba U$.  The
collection of all sinks (vertices without successors) is the {\it trivial}
coreset; the other coresets are {\it nontrivial}.
\sp

The definition immediately implies that no coreset contains another.
In Fig.~1, the set $\{x_7\}$ is the trivial coreset, and the nontrivial
coresets are $\{x_1,x_2,x_3,x_6\}$, $\{x_4\}$, and $\{x_5\}$.

\gpic{
\expandafter\ifx\csname graph\endcsname\relax \csname newbox\endcsname\graph\fi
\expandafter\ifx\csname graphtemp\endcsname\relax \csname newdimen\endcsname\graphtemp\fi
\setbox\graph=\vtop{\vskip 0pt\hbox{%
    \special{pn 8}%
    \special{ar 228 293 88 88 0 6.28319}%
    \special{ar 813 293 88 88 0 6.28319}%
    \special{ar 813 878 88 88 0 6.28319}%
    \special{ar 228 1463 88 88 0 6.28319}%
    \special{ar 1398 1463 88 88 0 6.28319}%
    \special{ar 1398 293 88 88 0 6.28319}%
    \special{ar 1983 1463 88 88 0 6.28319}%
    \graphtemp=.5ex\advance\graphtemp by 0.293in
    \rlap{\kern 0.228in\lower\graphtemp\hbox to 0pt{\hss $x_1$\hss}}%
    \graphtemp=.5ex\advance\graphtemp by 0.293in
    \rlap{\kern 0.813in\lower\graphtemp\hbox to 0pt{\hss $x_2$\hss}}%
    \graphtemp=.5ex\advance\graphtemp by 0.878in
    \rlap{\kern 0.813in\lower\graphtemp\hbox to 0pt{\hss $x_3$\hss}}%
    \graphtemp=.5ex\advance\graphtemp by 1.463in
    \rlap{\kern 0.228in\lower\graphtemp\hbox to 0pt{\hss $x_4$\hss}}%
    \graphtemp=.5ex\advance\graphtemp by 1.463in
    \rlap{\kern 1.398in\lower\graphtemp\hbox to 0pt{\hss $x_5$\hss}}%
    \graphtemp=.5ex\advance\graphtemp by 0.293in
    \rlap{\kern 1.398in\lower\graphtemp\hbox to 0pt{\hss $x_6$\hss}}%
    \graphtemp=.5ex\advance\graphtemp by 1.463in
    \rlap{\kern 1.983in\lower\graphtemp\hbox to 0pt{\hss $x_7$\hss}}%
    \special{pa 228 380}%
    \special{pa 228 977}%
    \special{fp}%
    \special{sh 1.000}%
    \special{pa 250 889}%
    \special{pa 228 977}%
    \special{pa 206 889}%
    \special{pa 250 889}%
    \special{fp}%
    \special{pa 228 977}%
    \special{pa 228 1375}%
    \special{fp}%
    \special{pa 290 355}%
    \special{pa 566 355}%
    \special{fp}%
    \special{sh 1.000}%
    \special{pa 479 333}%
    \special{pa 566 355}%
    \special{pa 479 377}%
    \special{pa 479 333}%
    \special{fp}%
    \special{pa 566 355}%
    \special{pa 751 355}%
    \special{fp}%
    \special{pa 751 230}%
    \special{pa 474 230}%
    \special{fp}%
    \special{sh 1.000}%
    \special{pa 562 252}%
    \special{pa 474 230}%
    \special{pa 562 209}%
    \special{pa 562 252}%
    \special{fp}%
    \special{pa 474 230}%
    \special{pa 290 230}%
    \special{fp}%
    \special{pa 813 790}%
    \special{pa 813 544}%
    \special{fp}%
    \special{sh 1.000}%
    \special{pa 791 632}%
    \special{pa 813 544}%
    \special{pa 835 632}%
    \special{pa 791 632}%
    \special{fp}%
    \special{pa 813 544}%
    \special{pa 813 380}%
    \special{fp}%
    \special{pa 315 1463}%
    \special{pa 912 1463}%
    \special{fp}%
    \special{sh 1.000}%
    \special{pa 824 1441}%
    \special{pa 912 1463}%
    \special{pa 824 1485}%
    \special{pa 824 1441}%
    \special{fp}%
    \special{pa 912 1463}%
    \special{pa 1310 1463}%
    \special{fp}%
    \special{pa 1486 1463}%
    \special{pa 1731 1463}%
    \special{fp}%
    \special{sh 1.000}%
    \special{pa 1644 1441}%
    \special{pa 1731 1463}%
    \special{pa 1644 1485}%
    \special{pa 1644 1441}%
    \special{fp}%
    \special{pa 1731 1463}%
    \special{pa 1895 1463}%
    \special{fp}%
    \special{pa 1398 1375}%
    \special{pa 1398 778}%
    \special{fp}%
    \special{sh 1.000}%
    \special{pa 1376 866}%
    \special{pa 1398 778}%
    \special{pa 1420 866}%
    \special{pa 1376 866}%
    \special{fp}%
    \special{pa 1398 778}%
    \special{pa 1398 380}%
    \special{fp}%
    \special{pa 1310 293}%
    \special{pa 1064 293}%
    \special{fp}%
    \special{sh 1.000}%
    \special{pa 1152 315}%
    \special{pa 1064 293}%
    \special{pa 1152 271}%
    \special{pa 1152 315}%
    \special{fp}%
    \special{pa 1064 293}%
    \special{pa 901 293}%
    \special{fp}%
    \special{pa 1398 205}%
    \special{pa 1105 0}%
    \special{pa 813 0}%
    \special{pa 520 0}%
    \special{pa 228 205}%
    \special{sp}%
    \special{pa 871 0}%
    \special{pa 783 0}%
    \special{fp}%
    \special{sh 1.000}%
    \special{pa 871 22}%
    \special{pa 783 0}%
    \special{pa 871 -21}%
    \special{pa 871 22}%
    \special{fp}%
    \special{pa 783 0}%
    \special{pa 725 0}%
    \special{fp}%
    \special{pa 166 355}%
    \special{pa 111 424}%
    \special{pa 19 384}%
    \special{pa -5 293}%
    \special{pa 19 201}%
    \special{pa 111 161}%
    \special{pa 166 230}%
    \special{sp}%
    \special{pa 13 340}%
    \special{pa 13 266}%
    \special{fp}%
    \special{sh 1.000}%
    \special{pa -8 353}%
    \special{pa 13 266}%
    \special{pa 35 353}%
    \special{pa -8 353}%
    \special{fp}%
    \special{pa 13 266}%
    \special{pa 13 216}%
    \special{fp}%
    \special{pa 751 940}%
    \special{pa 696 1009}%
    \special{pa 604 969}%
    \special{pa 579 878}%
    \special{pa 604 786}%
    \special{pa 696 746}%
    \special{pa 751 816}%
    \special{sp}%
    \special{pa 599 925}%
    \special{pa 599 851}%
    \special{fp}%
    \special{sh 1.000}%
    \special{pa 577 938}%
    \special{pa 599 851}%
    \special{pa 620 938}%
    \special{pa 577 938}%
    \special{fp}%
    \special{pa 599 851}%
    \special{pa 599 801}%
    \special{fp}%
    \hbox{\vrule depth1.580in width0pt height 0pt}%
    \kern 2.100in
  }%
}%
}
\ce{Fig.~1.  A digraph to illustrate coresets.}
\sp

Beineke and Zamfirescu [1] also introduced intersection digraphs, under the
name ``connection digraph''.  As studied in [8], a digraph $D$ is the {\it
intersection digraph} of a collection of pairs of sets
$\SET (S_v,T_v):{v\in V(D)}$ if $E(D)$ is the set
$\SET uv:{S_u\cap T_v\ne \nul}$.  Informally, we call $S_v$ and $T_v$ the {\it
source set} and {\it sink set} of $v$, and the edge $uv$ occurs when the
source set of $u$ intersects the sink set of $v$.  An alternative model
was introduced in [3].  Various classes of intersection digraphs 
have been studied.  Line digraphs are themselves intersection digraphs, arising
by letting $(S_e,T_e)=(\{y\},\{x\})$ for each edge $e=xy$.  We study a special
class of intersection digraphs that, like line digraphs, arise from arbitrary
digraphs and capture some of their structure.  Given a digraph $D$, the vertices
of this new digraph will be the coresets of $D$.

\SH
{2. PROPERTIES OF CORESETS}

\LM 1.
The coresets of a digraph are pairwise disjoint.
\PF
Since a sink is not a predecessor of any successor of itself, the trivial
coreset intersects no others.  Let $U,W$ be distinct nontrivial coresets.  Let
$X=U\cap W$.  Since $U$ is a coreset, $\ba X$ contains no element of $W-U$.
Since $W$ is a coreset, $\ba X$ contains no element of $U-W$.  Hence
$\ba X = X$.  Since $U$ and $W$ are minimal nonempty sets unchanged by
$\beta\comp\al$, we have $X=\nul$, and $U,W$ are disjoint.  \qed

The disjointness of distinct coresets implies that there is no edge from a
coreset to a successor of another coreset.  Equivalently, no two vertices in
distinct coresets have a common successor.

For convenience, in a digraph we define $\al(\nul)$ to be the collection of
source vertices and $\beta(\nul)$ to be the collection of sink vertices.
In describing partitions of a vertex set, we allow sets of the partition to
be empty.

\TH 2.
For every digraph $D$, the coresets of $D$ (including $\nul$) partition $V(D)$.
The successor sets of the coresets (including the set $\al(\nul)$ of sources)
also partition $V(D)$.
\PF
The first statement follows from Lemma 1 and the claim that every vertex $v$
belongs to a coreset.  We may assume that $v$ is not a sink.  Let $U$ be a
minimal set containing $v$ such that $\ba U =U$.  Such a set exists, since the
collection of all non-sinks has this property.  By the definition, $U$ contains
a coreset $X$.  If $v\notin X$, then $v$ cannot have a successor in $\al(X)$,
since $X$ is a coreset, and $\al(X)$ has no predecessor outside $X$.
Thus $\ba{U-X}=U-X$, which contradicts the minimality of $U$.  We conclude
that $v$ belongs to a coreset.

For the second statement, the observation that vertices of distinct coresets
cannot have common successors implies that the successor sets of the coresets
(and $\al(\nul)$) are pairwise disjoint.  Furthermore, every non-source vertex
is a successor of a vertex in some coreset, since the coresets partition $V(D)$.
\qed

In Fig.~1, the successor sets of the coresets listed earlier are
$\nul$, $\{x_1,x_2,x_3,x_4\}$, $\{x_5\}$, $\{x_6,x_7\}$, respectively.
This digraph has no source vertices.  In general, we would list
$U_0=\nul$ in the partition into coresets in order to obtain
the set $\al(\nul)$ of source vertices in the successor partition.

It is worth noting that the partitions in Theorem 2 are unchanged under reversal
of all edges.  Since $\ba U = U$ for a coreset $U$, the successor set $W=\al(U)$
of a coreset $U$ satisfies $\al(\beta(W))=W$.  Thus the coresets of the reversal
of $D$ are the successor sets in $D$ of the coresets in $D$.

Given a digraph $D$ with $A,B\esub V(D)$, we use $[A,B]$ to denote the set of
edges of $D$ from $A$ to $B$.  A {\it decomposition} of a digraph $D$ is a set
of pairwise edge-disjoint subgraphs whose union is $D$.

\TH 3.
Given a digraph $D$, let $\VEC U0m$ be the partition of $V(D)$ into
coresets, including $U_0=\nul$.  Let $D_i$ be the subdigraph of $D$
with vertex set $U_i\cup\al(U_i)$ and edge set $[U_i,\al(U_i)]$.
The digraphs $\VEC D0m$ decompose $D$.
\PF
By Theorem 2, every vertex and every edge appears in some $D_i$.
If some edge $uv$ appears in both $D_i$ and $D_j$ for $i\ne j$, then
$\{u,v\} \esub V(D_i)\cap V(D_j)$.  By Theorem 2, each of $\{u,v\}$
lies in $U_i\cap \al(U_j)$ or in $\al(U_i)\cap U_j$.  Since no coreset has an
edge to a successor of another, the two vertices cannot both lie in on of
these sets.  By symmetry, we may thus assume that $u\in U_i\cap \al(U_j)$ and
$v\in \al(U_i)\cap U_j$, which requires $u\ne v$.  Since $u\notin U_j$,
we now have $uv\notin E(D_j)$, a contradiction.  \qed

We call the resulting decomposition the {\it core decomposition} of $D$; the
subgraphs $\VEC D0m$ are the {\it core subgraphs}.  Fig.~2 shows the core
decomposition of the digraph in Fig.~1.

\gpic{
\expandafter\ifx\csname graph\endcsname\relax \csname newbox\endcsname\graph\fi
\expandafter\ifx\csname graphtemp\endcsname\relax \csname newdimen\endcsname\graphtemp\fi
\setbox\graph=\vtop{\vskip 0pt\hbox{%
    \special{pn 8}%
    \special{ar 219 281 84 84 0 6.28319}%
    \special{ar 781 281 84 84 0 6.28319}%
    \special{ar 781 844 84 84 0 6.28319}%
    \special{ar 219 1406 84 84 0 6.28319}%
    \special{ar 1344 281 84 84 0 6.28319}%
    \special{ar 359 1547 84 84 0 6.28319}%
    \special{ar 1344 1547 84 84 0 6.28319}%
    \special{ar 2188 1266 84 84 0 6.28319}%
    \special{ar 1484 1406 84 84 0 6.28319}%
    \special{ar 1484 422 84 84 0 6.28319}%
    \special{ar 2047 1406 84 84 0 6.28319}%
    \graphtemp=.5ex\advance\graphtemp by 0.281in
    \rlap{\kern 0.219in\lower\graphtemp\hbox to 0pt{\hss $x_1$\hss}}%
    \graphtemp=.5ex\advance\graphtemp by 0.281in
    \rlap{\kern 0.781in\lower\graphtemp\hbox to 0pt{\hss $x_2$\hss}}%
    \graphtemp=.5ex\advance\graphtemp by 0.844in
    \rlap{\kern 0.781in\lower\graphtemp\hbox to 0pt{\hss $x_3$\hss}}%
    \graphtemp=.5ex\advance\graphtemp by 1.406in
    \rlap{\kern 0.219in\lower\graphtemp\hbox to 0pt{\hss $x_4$\hss}}%
    \graphtemp=.5ex\advance\graphtemp by 0.281in
    \rlap{\kern 1.344in\lower\graphtemp\hbox to 0pt{\hss $x_6$\hss}}%
    \graphtemp=.5ex\advance\graphtemp by 0.562in
    \rlap{\kern 0.500in\lower\graphtemp\hbox to 0pt{\hss $D_1$\hss}}%
    \graphtemp=.5ex\advance\graphtemp by 1.547in
    \rlap{\kern 0.359in\lower\graphtemp\hbox to 0pt{\hss $x_4$\hss}}%
    \graphtemp=.5ex\advance\graphtemp by 1.547in
    \rlap{\kern 1.344in\lower\graphtemp\hbox to 0pt{\hss $x_5$\hss}}%
    \graphtemp=.5ex\advance\graphtemp by 1.266in
    \rlap{\kern 2.188in\lower\graphtemp\hbox to 0pt{\hss $x_7$\hss}}%
    \graphtemp=.5ex\advance\graphtemp by 1.406in
    \rlap{\kern 1.484in\lower\graphtemp\hbox to 0pt{\hss $x_5$\hss}}%
    \graphtemp=.5ex\advance\graphtemp by 0.422in
    \rlap{\kern 1.484in\lower\graphtemp\hbox to 0pt{\hss $x_6$\hss}}%
    \graphtemp=.5ex\advance\graphtemp by 1.406in
    \rlap{\kern 2.047in\lower\graphtemp\hbox to 0pt{\hss $x_7$\hss}}%
    \graphtemp=.5ex\advance\graphtemp by 1.458in
    \rlap{\kern 1.001in\lower\graphtemp\hbox to 0pt{\hss $D_2$\hss}}%
    \graphtemp=.5ex\advance\graphtemp by 0.984in
    \rlap{\kern 1.597in\lower\graphtemp\hbox to 0pt{\hss $D_3$\hss}}%
    \graphtemp=.5ex\advance\graphtemp by 1.125in
    \rlap{\kern 2.188in\lower\graphtemp\hbox to 0pt{\hss $D_4$\hss}}%
    \special{pa 219 366}%
    \special{pa 219 939}%
    \special{fp}%
    \special{sh 1.000}%
    \special{pa 240 855}%
    \special{pa 219 939}%
    \special{pa 198 855}%
    \special{pa 240 855}%
    \special{fp}%
    \special{pa 219 939}%
    \special{pa 219 1322}%
    \special{fp}%
    \special{pa 278 341}%
    \special{pa 544 341}%
    \special{fp}%
    \special{sh 1.000}%
    \special{pa 460 320}%
    \special{pa 544 341}%
    \special{pa 460 362}%
    \special{pa 460 320}%
    \special{fp}%
    \special{pa 544 341}%
    \special{pa 722 341}%
    \special{fp}%
    \special{pa 722 222}%
    \special{pa 456 222}%
    \special{fp}%
    \special{sh 1.000}%
    \special{pa 540 243}%
    \special{pa 456 222}%
    \special{pa 540 200}%
    \special{pa 540 243}%
    \special{fp}%
    \special{pa 456 222}%
    \special{pa 278 222}%
    \special{fp}%
    \special{pa 781 759}%
    \special{pa 781 523}%
    \special{fp}%
    \special{sh 1.000}%
    \special{pa 760 607}%
    \special{pa 781 523}%
    \special{pa 802 607}%
    \special{pa 760 607}%
    \special{fp}%
    \special{pa 781 523}%
    \special{pa 781 366}%
    \special{fp}%
    \special{pa 444 1547}%
    \special{pa 933 1547}%
    \special{fp}%
    \special{sh 1.000}%
    \special{pa 849 1526}%
    \special{pa 933 1547}%
    \special{pa 849 1568}%
    \special{pa 849 1526}%
    \special{fp}%
    \special{pa 933 1547}%
    \special{pa 1259 1547}%
    \special{fp}%
    \special{pa 1569 1406}%
    \special{pa 1805 1406}%
    \special{fp}%
    \special{sh 1.000}%
    \special{pa 1721 1385}%
    \special{pa 1805 1406}%
    \special{pa 1721 1427}%
    \special{pa 1721 1385}%
    \special{fp}%
    \special{pa 1805 1406}%
    \special{pa 1963 1406}%
    \special{fp}%
    \special{pa 1484 1322}%
    \special{pa 1484 832}%
    \special{fp}%
    \special{sh 1.000}%
    \special{pa 1463 917}%
    \special{pa 1484 832}%
    \special{pa 1505 917}%
    \special{pa 1463 917}%
    \special{fp}%
    \special{pa 1484 832}%
    \special{pa 1484 506}%
    \special{fp}%
    \special{pa 1259 281}%
    \special{pa 1023 281}%
    \special{fp}%
    \special{sh 1.000}%
    \special{pa 1108 302}%
    \special{pa 1023 281}%
    \special{pa 1108 260}%
    \special{pa 1108 302}%
    \special{fp}%
    \special{pa 1023 281}%
    \special{pa 866 281}%
    \special{fp}%
    \special{pa 1344 197}%
    \special{pa 1063 0}%
    \special{pa 781 0}%
    \special{pa 500 0}%
    \special{pa 219 197}%
    \special{sp}%
    \special{pa 838 0}%
    \special{pa 753 0}%
    \special{fp}%
    \special{sh 1.000}%
    \special{pa 838 21}%
    \special{pa 753 0}%
    \special{pa 838 -20}%
    \special{pa 838 21}%
    \special{fp}%
    \special{pa 753 0}%
    \special{pa 697 0}%
    \special{fp}%
    \special{pa 159 341}%
    \special{pa 106 408}%
    \special{pa 19 369}%
    \special{pa -5 281}%
    \special{pa 19 193}%
    \special{pa 106 155}%
    \special{pa 159 222}%
    \special{sp}%
    \special{pa 13 327}%
    \special{pa 13 255}%
    \special{fp}%
    \special{sh 1.000}%
    \special{pa -7 340}%
    \special{pa 13 255}%
    \special{pa 34 340}%
    \special{pa -7 340}%
    \special{fp}%
    \special{pa 13 255}%
    \special{pa 13 208}%
    \special{fp}%
    \special{pa 722 903}%
    \special{pa 669 970}%
    \special{pa 581 932}%
    \special{pa 556 844}%
    \special{pa 581 756}%
    \special{pa 669 717}%
    \special{pa 722 784}%
    \special{sp}%
    \special{pa 575 889}%
    \special{pa 575 818}%
    \special{fp}%
    \special{sh 1.000}%
    \special{pa 554 902}%
    \special{pa 575 818}%
    \special{pa 596 902}%
    \special{pa 554 902}%
    \special{fp}%
    \special{pa 575 818}%
    \special{pa 575 770}%
    \special{fp}%
    \hbox{\vrule depth1.659in width0pt height 0pt}%
    \kern 2.300in
  }%
}%
}
\ce{Fig.~2.  A core decomposition.}
\sp

\PR 4.
If $D'$ is a core subgraph of a digraph $D$, then the core decomposition of
$D'$ has only one nonempty subgraph, $D'$ itself.
\PF
Let $U$ be the coreset of $D$ such that $D'$ consists of all edges in $D$
from $U$ to $\al(U)$.  We may assume that $U$ is a nontrivial coreset of $D$,
so that $\al(U)-U$ is the set of sinks in $D'$.  Suppose $U'$ is a nontrivial
coreset in $D'$, so $\beta_{D'}(\al_{D'}(U'))=U'$.  Because $\al(U)$ contains
all successors in $D$ of vertices of $U'$, we have $\al_{D'}(U')=\al(U')$.
Similarly, since $U$ is a coreset of $D$, $U$ contains all predecessors of
successors of vertices of $U'$, so $\beta_{D'}(\al_{D'}(U'))=\ba{U'}$.
We have proved that $\ba{U'}=U'$.  By the minimality property of coresets,
$U'=U$.  Thus the coresets of $D'$ are $U$ and $\al(U)-U$, and the only
core subgraph is $D'$ itself.  \qed

We can describe coresets in the language of adjacency matrices.

\PR 5.
A nonempty vertex set $U$ in a digraph $D$ is a nontrivial coreset if and only
if the rows of the adjacency matrix that correspond to $U$ are nonzero, are
orthogonal to the remaining rows, and have no nonempty subset with the same
property.
\PF
The conditions state that the vertices of $U$ are not sinks, that 
$\ba U\esub U$, and that $U$ is a minimal nonempty set with these properties.
Since $u\in\ba u$ whenever $u$ is not a sink, these are precisely the
conditions for a nontrivial coreset.  \qed

In the language of adjacency matrices, Theorem 3 states that if the
rows of $D$ are grouped as $\VEC U0m$ and the columns as 
$\al(U_0),\ldots,\al(U_m)$, then the adjacency matrix $A(D)$ partitions into
blocks $B_{i,j}$ with $0\le i,j\le m$ such that $B_{i,i} = A(D_i)$ and all
other blocks are zero matrices.

\SH
{3. CHARACTERIZATION OF LINE DIGRAPHS}

Geller and Harary [2] proved that a digraph $D$ is a line digraph if and 
only if there exist partitions $\VEC A1r$ and $\VEC B1r$ of $V(D)$ (using
possibly empty sets) such that $E(D) = \UE i1r A_i\times B_i$.  We give several
related characterizations using the concept of coreset, including the
Geller-Harary characterization for completeness.  A digraph with vertex
set $A\cup B$ and edge set $A\times B$ is a {\it Cartesian product digraph}.

\TH 6.
For a digraph $D$ with core decomposition $\VEC D0m$, the following are
equivalent:
\br A) $D$ is a line digraph.
\br B) $V(D)$ admits partitions $\{A_i\}$ and $\{B_i\}$ such that
$E(D) = \cup (A_i\times B_i)$.
\br C) $\al(u)=\al(x)$ whenever $u$ and $x$ belong to the same coreset in $D$.
\br D) For each $uv\in E(D)$, the set $\beta(v)$ is the coreset containing $u$.
\br E) Each $D_i$ is a Cartesian product digraph.
\br F) Each $D_i$ is a line digraph.
\PF
A$\Imp$B (Geller-Harary [2]): Given $D=L(D')$, where $V(D')=\{\VEC w1r\}$, let
$A_i$ be the set of edges in $D'$ with head $w_i$, and let $B_i$ consist of
those with tail $w_i$ (when $w_i$ is a source or sink, $A_i$ or $B_i$ is empty,
respectively).  These sets also partition $V(D)$, and the definition of line
digraph yields $E(D) = \UE i1r A_i\times B_i$.
 
B$\Imp$C: Given the partitions as described, a vertex of $A_i$ can have no
successors outside $B_i$, and a vertex of $B_i$ can have no predecessors outside
$A_i$.  Thus $A_i$ is a coreset and $B_i$ is its successor set, and $\al(x)=B_i$
for all $x\in A_i$.

C$\Imp$D: Given an edge $uv$, Condition C implies that $\beta(v)$ contains the
coreset that contains $u$; by the definition of coreset, it cannot contain more.

D$\Imp$E: Condition D implies that each vertex of a coreset $U$ is a predecessor
of each vertex of the successor set $\al(U)$, and thus
$E(D_i)=U_i\times\al(U_i)$.

E$\Imp$A: As a Cartesian product, we must have $E(D_i)=U_i\times\al(U_i)$.
We must assume that the sink set is $\al(U_m)$ and the set $U_0$ is empty, so
that $\al(U_0)$ is the source set.  Construct a digraph $D'$ with
$V(D')=\VEC w0m$.  Put an edge in $D'$ for each $v\in V(D)$.  The edge $v$ in
$D'$ is $w_iw_j$ if $v$ belongs to $U_j$ in the coreset partition and to
$\al(U_i)$ in the successor partition ($D'$ may have multiple edges).
By construction, $L(D')=D$.

E$\Imp$F$\Imp$D: Every Cartesian product digraph is explicitly a line digraph.
Finally, F$\Imp$D follows by applying Proposition 4 and the equivalence of A and
D to each $D_i$.
\qed

The expression of Theorem 3 in the language of adjacency matrices allows us
to interpret Theorem 6 in these terms also.  In particular, a digraph is
a line digraph if and only if its rows and columns can be permuted so that
the 1's form rectangular blocks that do not share rows or columns.  These
immediately implies the characterization by Richards [7] that a 0,1-matrix
is the adjacency matrix of a line digraph if and only if any two columns
(or any two rows) are identical or orthogonal.  We could also translate
conditions C and D of Theorem 6 into conditions on the rows corresponding
to vertices in a given coreset.

\SH
{4. $n$th-ORDER LINE DIGRAPHS}
Introduced in [1],
the {\it $n$-order line digraph} $L^n(D)$ of a digraph $D$ is the digraph
obtained from $D$ by iteratively applying the line digraph operator $n$ times.
Coresets enable us to characterize the $n$th-order line digraphs of digraphs
with no sources or sinks.  Note that the vertex set of $L^n(D)$ is the set of
$n$-walks (walks of length $n$) in $D$, with an edge from an $n$-walk $W_1$ to
an $n$-walk $W_2$ if deleting the first vertex and edge from $W_1$ yields the
same $n-1$-walk as deleting the last vertex and edge from $W_2$.  We confine our
attention to digraphs without sources or sinks to avoid annoying technicalities
about the nonexistence of walks of given lengths from a given vertex.

For the characterization, we introduce $n$-th order coresets.  Given a digraph
$D$, define the digraph $D^n$ by $V(D^n)=V(D)$ and $E(D^n)=\{ uv\st{D}$ has a
walk of length $n$ from $u$ to $v\}$.  The successor operator in $D^n$ is the
$n$th iterate of the successor operator in $D$: $\al_{D^n}(u)=\al^n(u)$.
Similarly $\beta_{D^n}(u)=\beta^n(u)$.  The {\it $n$th-order coresets} of $D$
are the coresets of $D^n$.  Equivalently, the $n$th-order coresets of a digraph
without sources or sinks are the minimal nonempty sets $U$ such that
$\beta^n(\al^n(U))=U$.  By applying Theorem 2 to $D^n$, we
immediately obtain the following corollary.
\skipit{
The final statement of the
corollary can also be phrased as saying that if $u$ and $v$ belong to distinct
$n$th-order coresets, then there do not exist $n$-walks from both $u$ and $v$
to a common vertex $w$.
}

\CO 7.
If $D$ is a digraph without sources or sinks, then
the $n$th-order coresets partition $V(D)$, as do the
$n$th-order successor sets of these coresets.  Furthermore,
vertices of distinct $n$th-order coresets do not have a common
$n$th-order successor.  \qed\nosp

\LM 8.
If $D$ is a digraph without sources or sinks, then $U$ is an $n$th-order coreset
in $D$ if and only if $W$ is an $(n+1)$th-order coreset in $L(D)$, where $W$ is
the set of edges in $D$ with heads in $U$.
\PF
Observe first that an $m$th-order coreset $W$ in $L(D)$ contains all edges of
$D$ that share heads with its members.  If $uv\in W$ and $xv\in E(D)$, then
$xv$ can be substituted for $uv$ as the first edge in an $m+1$-walk in $D$.
Thus $xv\in \beta^m(\al^m(W))=W$ in $L(D)$.  As a result, if $W$ is an 
$(n+1)$th-order coreset in $L(D)$ and $U$ is the subset of $V(D)$ consisting
of heads of elements of $W$, then $W$ is obtained from $U$ as defined above.

We prove that $U$ is a coreset in $D^n$ if and only if $W$ is a coreset in
$[L(D)]^{n+1}$.  The digraph $[L(D)]^{n+1}$ has an edge from $uv\in E(D)$ to
$xy\in E(D)$ if and only if $L(D)$ has a walk of length $n+1$ from $uv$ to
$xy$.  Such a walk exists if and only if $D$ has an $n$-walk from $v$ to $x$,
which corresponds to existence of the edge $vx$ in $D^n$.

Thus $z$ is a predecessor of a successor of $v$ in $D^n$ if and only if
there exist edges $uv$ and $tz$ in $E(D)$ such that $tz$ is a predecessor of
a successor of $uv$ in $[L(D)]^{n+1}$.  This implies that$\ba U=U$ in $D^n$
if and only if $\ba W=W$ in $[L(D)]^{n+1}$.  The equality holds for a 
proper subset of $U$ if and only if it holds for a proper subset of $W$;
this completes the proof of the claim.  \qed\nosp

\TH 9.
If $D$ is a digraph without sources or sinks, then $D$ is an $n$th-order
line digraph if and only if, for every $1\le i\le n$ and every $i$th-order
coreset $U$, there is exactly one $i$-walk from each vertex of $U$ to
each vertex of $\al^i(U)$.
\PF
We call the desired condition for $i$ the {\it $i$-uniqueness condition} in $D$.
We first compare $m$-uniqueness in $D$ with $(m+1)$-uniqueness in $L(D)$.  As
shown in Lemma 8, for each $m$th-order coreset $U$ in $D$ there is an
$(m+1)$th-order coreset $W$ in $L(D)$ (and vice versa) such that each vertex of
$W$ is an edge of $D$ whose head is in $U$.  Thus also each vertex of
$\al_{L(D)}^{m+1}(W)$ is an edge of $D$ whose tail is in $\al_D^m(U)$.   Thus
$(m+1)$-walks in $L(D)$ beginning in $W$ correspond naturally to $m$-walks in
$D$ beginning in $U$.  In particular, the $m$-uniqueness condition in $D$
is equivalent to the $(m+1)$-uniqueness condition in $L(D)$.

The result now follows easily by induction on $n$.  The case $n=1$ is contained
in Theorem 6.  For the induction step, suppose that the characterization holds
when $n=m$; we prove it for $n=m+1$.  For necessity, suppose that $D$ is an
$m$th-order line digraph, so $L(D)$ is an $(m+1)$th-order line graph.
By the induction hypothesis, the $i$-uniqueness condition holds in $D$ for
$1\le i\le m$, and hence the $i$-uniqueness condition holds in $L(D)$ for
$2\le i\le m+1$.  Since $L(D)$ is a line digraph, it also holds in $L(D)$ for
$i=1$.

Conversely, suppose that $D$ satisfies the $i$-uniqueness condition for
$1\le i\le m+1$.  Since the claim holds for $n=1$, we have $D=L(D')$ for some
digraph $D'$.  Now $D'$ satisfies the $i$-uniqueness condition for 
$1\le i\le m$.  By the induction hypothesis, $D'$ is an $m$th-order line
digraph, and thus $D$ is an $(m+1)$th-order line digraph.  \qed

\SH
{5. THE CORESET DIGRAPH OF A DIGRAPH}

By Theorem 2, the coresets of a digraph partition its vertex set.  We define
a special derived digraph of $D$ with the elements of the coreset partition as
the vertex set.  Like the line digraph $D$, it is an intersection digraph where
the sets used in the intersection representation are subsets of $V(D)$. 

\DF{}
If $D$ is a digraph with coreset partition $\VEC U0m$, then the 
{\it coreset digraph} $Y(D)$ is the digraph without multiple edges defined by
$V(Y(D))= \VEC U0m$ and $E(Y(D))=\SET U_iU_j:{\al(U_i)\cap U_j\ne\nul}$.
We iterate this operation by defining $Y^1(D)=Y(D)$ and $Y^n(D)=Y(Y^{n-1}(D))$
for $n\ge2$.
\sp

Fig.~3 shows the coreset digraph of the digraph in Fig.~1.
Note that $Y(D)$ is the intersection digraph arising when the pair
$(S_i,T_i)$ assigned to the vertex $U_i$ is $(\al(U_i),U_i)$.
We give another interpretation of this operation and describe
the effect of iterating it.

\gpic{
\expandafter\ifx\csname graph\endcsname\relax \csname newbox\endcsname\graph\fi
\expandafter\ifx\csname graphtemp\endcsname\relax \csname newdimen\endcsname\graphtemp\fi
\setbox\graph=\vtop{\vskip 0pt\hbox{%
    \special{pn 8}%
    \special{ar 247 128 95 95 0 6.28319}%
    \special{ar 564 1079 95 95 0 6.28319}%
    \special{ar 1198 762 95 95 0 6.28319}%
    \special{ar 1673 762 95 95 0 6.28319}%
    \graphtemp=.5ex\advance\graphtemp by 0.128in
    \rlap{\kern 0.247in\lower\graphtemp\hbox to 0pt{\hss $U_1$\hss}}%
    \graphtemp=.5ex\advance\graphtemp by 1.079in
    \rlap{\kern 0.564in\lower\graphtemp\hbox to 0pt{\hss $U_2$\hss}}%
    \graphtemp=.5ex\advance\graphtemp by 0.762in
    \rlap{\kern 1.198in\lower\graphtemp\hbox to 0pt{\hss $U_3$\hss}}%
    \graphtemp=.5ex\advance\graphtemp by 0.762in
    \rlap{\kern 1.673in\lower\graphtemp\hbox to 0pt{\hss $U_4$\hss}}%
    \special{pa 247 223}%
    \special{pa 437 679}%
    \special{fp}%
    \special{sh 1.000}%
    \special{pa 422 582}%
    \special{pa 437 679}%
    \special{pa 378 601}%
    \special{pa 422 582}%
    \special{fp}%
    \special{pa 437 679}%
    \special{pa 564 984}%
    \special{fp}%
    \special{pa 631 1012}%
    \special{pa 931 902}%
    \special{fp}%
    \special{sh 1.000}%
    \special{pa 833 912}%
    \special{pa 931 902}%
    \special{pa 849 957}%
    \special{pa 833 912}%
    \special{fp}%
    \special{pa 931 902}%
    \special{pa 1130 829}%
    \special{fp}%
    \special{pa 1130 695}%
    \special{pa 641 395}%
    \special{fp}%
    \special{sh 1.000}%
    \special{pa 709 465}%
    \special{pa 641 395}%
    \special{pa 734 424}%
    \special{pa 709 465}%
    \special{fp}%
    \special{pa 641 395}%
    \special{pa 314 195}%
    \special{fp}%
    \special{pa 1293 762}%
    \special{pa 1464 762}%
    \special{fp}%
    \special{sh 1.000}%
    \special{pa 1369 738}%
    \special{pa 1464 762}%
    \special{pa 1369 786}%
    \special{pa 1369 738}%
    \special{fp}%
    \special{pa 1464 762}%
    \special{pa 1578 762}%
    \special{fp}%
    \special{pa 179 195}%
    \special{pa 120 270}%
    \special{pa 21 227}%
    \special{pa -6 128}%
    \special{pa 21 29}%
    \special{pa 120 -14}%
    \special{pa 179 61}%
    \special{sp}%
    \special{pa 15 179}%
    \special{pa 15 98}%
    \special{fp}%
    \special{sh 1.000}%
    \special{pa -8 194}%
    \special{pa 15 98}%
    \special{pa 38 194}%
    \special{pa -8 194}%
    \special{fp}%
    \special{pa 15 98}%
    \special{pa 15 45}%
    \special{fp}%
    \hbox{\vrule depth1.206in width0pt height 0pt}%
    \kern 1.800in
  }%
}%
}
\ce{Fig.~3.  A coreset digraph.}
\sp

\PR {10}.
For a digraph $D$, $Y(D)$ is the digraph obtained by identifying the vertices of
each coreset of $D$ into a single vertex and deleting the resulting extra copies
of edges.  Furthermore, $Y(D)=D$ if and only if each nonempty coreset of $D$ is
a single vertex, and the sequence $Y^n(D)$ always converges to some digraph $F$
as $n\to\infty$ (we write $Y^n(D)\to F$).
\PF
The coreset digraph of $D$ has a edge from $U_i$ to $U_j$ precisely when some
vertex of $U_i$ has a successor in $U_j$; this is achieved by identifying
vertices within coresets.  Thus the order of $Y(D)$ is less than that of
$D$ if and only some coreset has at least two elements, and the characterization
of $Y(D)=D$ follows.  Since the order of a digraph is an integer, the sequence
$Y^n(D)$ converges because $Y(D)\ne D$ if and only if $Y(D)$ has fewer vertices
than $D$.  \qed\upsp

\CO {11}.
A digraph $D$ is isomorphic to its coreset digraph if and only if it has at most
one sink and has maximum indegree at most one.  In particular, $Y(D)=D$ if
and only if $D$ is a path or is a digraph obtained by identifying a vertex of
a cycle with the source of a path.  Furthermore, $Y^n(D)$ always converges
to such a graph.
\PF
If there is more than one sink or some vertex has more than one predecessor,
then some coreset has size at least two.  Conversely, if the condition holds,
then each coreset has size one.  The digraphs described are the only ones 
where the condition holds.  Finally, $Y^n(D)$ always converges to a digraph $F$
such that $Y(F)=F$.  \qed

One could measure the complexity of $D$ by the value $n$ where $Y^n(D)$ first
reaches its limit.

\SH
{\ce{Acknowledgments}}
The first author is very grateful to Professor Jinsheng He for his encouragement
and help.  He is also grateful to Miss Zhui for her help.

\SH
{\ce{References}}
\frenchspacing
\BP[1]
L.W.~Beineke and C.M.~Zamfirescu, Connection digraphs and second order line
graphs.  \DM\ 39(1982), 237--254.
\BP[2]
D.P. Geller and F. Harary, Arrow diagrams are line digraphs,
\SIAP\ 16(1968), 1141--1145.
\BP[3]
F.~Harary, J.A.~Kabell, and F.R.~McMorris, Bipartite intersection graphs,
{\it Comm.\ Math.\ Univ.\ Carolinae} 23(1982), 739--745.
\BP[4]
F. Harary and R.Z. Norman, Some properties of line digraphs,
{\it Rend. Circ. Mat. Palermo} 9(1960), 161--168.
\BP[5]
R.L. Hemminger and L.W. Beineke, Line graphs and line digraphs,
in {\it Selected Topics in Graph Theory} (L.W. Beineke and R.J. Wilson, eds.).
Academic Press (1978), 271--305.
\BP[6]
C. Heuchenne, Sur une certaine correspondance entre graphes,
{\it Bull. Soc. Roy. Sci. Liege} 33(1964), 743--753.
\BP[7]
P.I. Richards, Precedence constraints and arrow diagrams,
\SIREV\ 9(1967), 548--553.
\BP[8]
M. Sen, S. Das, A.B. Roy, and D.B. West, Interval digraphs: an analogue of
interval graphs, \JGT\ 13(1989), 189-202.
\bye